\documentstyle[12pt]{article}

\newtheorem{defin}{Definition}[section]

\newtheorem{prop}[defin]{Proposition}
\newtheorem{theorem}[defin]{Theorem}
\newtheorem{cor}[defin]{Corollary}
\newtheorem{conj}[defin]{Conjecture}

\newcommand{\be}{\begin{equation}}
\newcommand{\ee}{\end{equation}}
\newcommand{\ba}{\begin{eqnarray}}
\newcommand{\ea}{\end{eqnarray}}
\newcommand{\n}{\nonumber \\}

\newcommand{\eq}[1]{(\ref{#1})}

\newcommand{\ket}[1]{|#1\rangle}
\newcommand{\bra}[1]{\langle#1|}

\newcommand\bZ{\bf Z}

\def\generalYoung{
\vskip.25cm
\noindent
\makebox[  4cm]{ }
\makebox[  2cm]{$s_N$}\hskip-.4pt
\makebox[1.7cm]{$s_{N-1}$}
\makebox[1.4cm]{$\!\cdots\cdot$}
\makebox[1.3cm]{$\!\!s_2$}\hskip-.5pt
\makebox[1.3cm]{$\!\!s_1$}
\hfill\break
 \makebox[  4cm][r]{$\hfill\lambda\:=\;$}
\framebox[  2cm][l]{\rule[  -1cm]{0cm}{  2cm}
                    \raisebox{.2cm}{$\!r_N$}}\hskip-.4pt
\framebox[1.7cm][l]{\rule[-0.7cm]{0cm}{1.7cm}
                    \raisebox{.2cm}{$\!r_{N-1}$}}
\hskip-0.15cm\rule[1.105cm]{1.6cm}{0.4pt}\hskip-1.55cm
 \makebox[1.4cm]   {\raisebox{.2cm}{$\,\,\cdots\cdot$}}
\framebox[1.3cm][l]{\rule[-0.3cm]{0cm}{1.3cm}
                    \raisebox{.2cm}{$\!r_2$}}\hskip-.4pt
\framebox[1.3cm][l]{\rule[0.0cm]{0cm}{1.0cm}
                    \raisebox{.2cm}{$\!r_1$}}
\makebox[0.5cm][r]{,}
\vskip.3cm
}

\begin{document}
\renewcommand{\thefootnote}{\fnsymbol{footnote}}
\font\csc=cmcsc10 scaled\magstep1
{\baselineskip=14pt
 \rightline{
 \vbox{\hbox{YITP/U-95-19}
       \hbox{DPSU-95-1}
       \hbox{UT-706}
       \hbox{May 1995}
}}}

\vskip 5mm
\begin{center}
{\large\bf
Integral Representations\\ of\\ the Macdonald Symmetric Functions}

\vspace{10mm}

{\csc Hidetoshi AWATA}\footnote{JSPS fellow}\setcounter{footnote}
{0}\renewcommand{\thefootnote}{\arabic{footnote}}\footnote{
      e-mail address : awata@yisun1.yukawa.kyoto-u.ac.jp},
{\csc Satoru ODAKE}\footnote{
      e-mail address : odake@yukawa.kyoto-u.ac.jp}\\
\vskip.1in
and \
{\csc Jun'ichi SHIRAISHI}$^*$\footnote{
      e-mail address : shiraish@danjuro.phys.s.u-tokyo.ac.jp}

{\baselineskip=15pt
\it\vskip.25in
  $^{1}$Uji Research Center, Yukawa Institute for Theoretical Physics \\
  Kyoto University, Uji 611, Japan \\
\vskip.1in
  $^2$Department of Physics, Faculty of Science \\
  Shinshu University, Matsumoto 390, Japan \\
\vskip.1in
  $^3$Department of Physics, Faculty of Science \\
  University of Tokyo, Tokyo 113, Japan
}

\end{center}

\vspace{4mm}

\begin{abstract}
{
Multiple-integral representations of the (skew-)Macdonald
symmetric functions are obtained.
Some bosonization schemes for the integral representations are also
constructed.
}
\end{abstract}

\vskip 3cm
q-alg/9506006
\setcounter{footnote}{0}
\renewcommand{\thefootnote}{\arabic{footnote}}
\newpage
\section{Introduction}

Calogero--Sutherland model\cite{rCS} and its various
generalizations\cite{rHS,rKHH} have been extensively studied and
these $1/r^2$ type models are known to describe systems
with the generalized exclusion principle in $1+1$ dimension\cite{rHal}.
Calogero--Sutherland model describes a system of non--relativistic
particles on a circle under the inverse square potential.
Its Hamiltonian and momentum are
\be
  H_{CS}=\sum_{j=1}^{N_0}\frac{1}{2}
  \biggl(\frac{1}{i}\frac{\partial}{\partial q_j}\biggr)^2
  +\Bigl(\frac{\pi}{L}\Bigr)^2
  \sum_{i,j=1 \atop i<j}^{N_0}
  \frac{\beta(\beta-1)}{\sin^2\frac{\pi}{L}(q_i-q_j)},
  \quad
  P_{CS}=\sum_{j=1}^{N_0}
  \frac{1}{i}\frac{\partial}{\partial q_j}.
  \label{CS}
\ee
where $\beta$ is a coupling constant.
This Calogero--Sutherland model is related to many branches of
low--dimensional physics and mathematics;
quantum Hall effect\cite{rQHE},
2D Yang-Mills theory\cite{rGN,rMP1},
matrix model\cite{rAJ,rSLA},
Yangian symmetry\cite{rHHTBP,rBGHP},
Virasoro symmetry\cite{rAMOS1,rMY1,rMY2},
$W_{1+\infty}$ symmetry\cite{rUWH},
Laplace-Beltrami operator \cite{rHO,rCM}, etc.
One of the recent remarkable development was the evaluation of
some dynamical correlation functions\cite{rSLA,rF,rHa,rLPS,rMP2}.
In these calculation the Jack symmetric functions\cite{rMac,rS}
play a central role, because they describe the excited states
of the Calogero-Sutherland model.
In the previous works\cite{rAMOS1,rMY1,rMY2,rAMOS2}, the free field
realization of the wave functions, in other words, the integral representations
of the Jack symmetric functions is discussed.

Several years ago,
Ruijsenaars constructed a relativistic (or lattice regularized)
version of the Calogero system\cite{rRui}.
That model is integrable, since it
has mutually commuting hermitian operators
$\hat{S}_k$ ($k=1,2,\cdots,N_0$):
\begin{equation}
\hat{S}_k=\sum_{I\subset \{ 1,\cdots,N_0 \} \atop |I|=k}
    \prod_{i\in I \atop j \not\in I }h(q_j-q_i)^{\frac12}\cdot
    \exp \Biggl( \rho \sum_{j\in I}
           \frac{1}{i}\frac{\partial}{\partial q_j}\Biggr)
    \prod_{i\in I \atop j \not\in I }h(q_i-q_j)^{\frac12},
\end{equation}
where $h(q)=\sigma(q+\mu)/\sigma(q)$ and
$\rho \in {\bf R}$, $\mu\in{\bf C}$.
Here $\sigma (z)$ denotes the Weierstrass $\sigma$-function
defined by
$$
\sigma(z)=z \prod_{m,n\in{\bf Z}\atop (m,n)\neq(0,0)}
  \left(1-\frac{z}{\Omega_{m,n}}\right)
  \exp\left( \frac{z}{\Omega_{m,n}}+\frac{z^2}{2\Omega_{m,n}}\right),
$$
where $\Omega_{m,n}=2m\omega_1+2n\omega_3$ and $2\omega_1$ and
$2\omega_3$ denote primitive periods.
The relation to the Calogero-Sutherland model is the
following:
\begin{eqnarray}
\lim_{\rho \rightarrow 0}\frac{1}{\rho^2}\left(
\frac{\hat{S}_1+\hat{S}_{N_0}^{-1}\hat{S}_{N_0-1}}{2}-N_0\right)
&\!\!=\!\!&\sum_{j=1}^{N_0}\frac{1}{2}
  \biggl(\frac{1}{i}\frac{\partial}{\partial q_j}\biggr)^2
  +  \beta(\beta-1)\sum_{i,j=1 \atop i<j}^{N_0}
  \wp(q_i-q_j),\\
\lim_{\rho \rightarrow 0}\frac{1}{\rho}\left(
\frac{\hat{S}_1-\hat{S}_{N_0}^{-1}\hat{S}_{N_0-1}}{2}\right)
&\!\!=\!\!&\sum_{j=1}^{N_0}
  \frac{1}{i}\frac{\partial}{\partial q_j},
\end{eqnarray}
where we set $\mu=i\beta\rho$ and $\wp(z)$ is the Weierstrass
$\wp$-function given by
$\wp(z)=-\frac{d}{dz}\zeta(z)$,
$\zeta(z)=\frac{d}{dz}\sigma(z)/\sigma(z)$.
If we consider the case $2\omega_1=L$, $\omega_3=i\infty$,
we have
$$
\wp(z)=\Bigl(\frac{\pi}{L}\Bigr)^2
  \left(\frac{1}{\sin^2\frac{\pi}{L}z}-\frac{1}{3}\right).
$$
Thus the system reduces to the model defined by (\ref{CS}).

To solve Ruijsenaars' system, we need an explicit formula
for the simultaneous eigenfunctions of $\hat{S}_k$'s.
When the $\wp$-function degenerates to the trigonometric function,
the commuting operators $\hat{S}_k$'s essentially degenerate
Macdonald's operators\cite{rMac}.
Therefore, the eigenfunctions are given by
the Macdonald symmetric functions.
(As for definitions of the Macdonald symmetric functions,
see the following sections.)
In this article, we construct integral representations of
the Macdonald symmetric functions and construct some boson
realization schemes of the integral formula.
These results are considered as natural deformation theories
of the previous works on the Jack symmetric
functions\cite{rAMOS1,rMY1,rMY2,rAMOS2}.
We hope that the general elliptic case
may be treated in a similar manner.

This paper is organized as follows.
In section 2 we give a short summary of
the Macdonald symmetric functions.
In section 3  we derive multiple integral representation
formulas for the Macdonald symmetric functions by using
two types of maps.
Moreover using the isomorphism between the ring of symmetric functions
and the boson Fock space, we derive the integral representations
of the skew Macdonald functions and the Kostka matrix.
In section 4 we construct two bosonization formulas
for the integral representation of the Macdonald
symmetric functions obtained in the section 3.
Two cases ($\beta\in{\bf Z}_{>0}$ and $\beta\in{\bf C}$)
are discussed separately. $A$-type structure and
finite temperature calculation method are used respectively.
Section 5 is devoted to discussions.

\section{Brief Review of Macdonald Symmetric Functions}
\setcounter{equation}{0}

In this section we review some basic facts about
the Macdonald symmetric functions \cite{rMac}.

\noindent
{\large\bf 2.1}~{\sl Notations and the scalar product
$\langle\;,\;\rangle_{q,t}$}\\
Let $q$, $t$ be independent indeterminates and
$\Lambda_{n;{\bf Q}(q,t)}$ be the ring of
symmetric functions in $n$ variables ($x_i$) over
the field of rational functions in $q$ and $t$.
We sometimes write $t=q^{\beta}$.
In the limit of $q\rightarrow 1$ this $\beta$ is
understood as the coupling constant of the Calogero-Sutherland
model.
The ring of symmetric functions $\Lambda_{{\bf Q}(q,t)}$ is defined by
$$
  \Lambda_{{\bf Q}(q,t)}=
  \lim_{\scriptstyle \longleftarrow \atop \scriptstyle n}
  \Lambda_{n,{\bf Q}(q,t)},
$$
where the projective limit is given by the following homomorphism:
\begin{eqnarray*}
\Lambda_{n+1;{\bf Q}(q,t)}&\rightarrow& \Lambda_{n;{\bf Q}(q,t)},\\
f(x_1,\cdots,x_n,x_{n+1})&\mapsto& f(x_1,\cdots,x_n,0).
\end{eqnarray*}

There are various bases of the ring $\Lambda_{{\bf Q}(q,t)}$.
They are indexed by partitions.
A partition $\lambda$ is a sequence
$\lambda=(\lambda_1,\lambda_2,\lambda_3,\cdots)$
of non-negative integers, such that
$\lambda_1 \geq\lambda_2 \geq\cdots$ and
$|\lambda|=\sum_i \lambda_i < \infty$.
The nonzero $\lambda_i$'s are called the parts of $\lambda$,
and the number of parts is the length $\ell(\lambda)$ of $\lambda$.
For two partitions $\lambda,\mu$, we define
$\lambda+\mu=(\lambda_1+\mu_1,\lambda_2+\mu_2,\cdots)$.
The natural partial ordering  is defined as follows:
\be
\lambda\geq\mu \Leftrightarrow
|\lambda|=|\mu| \;{\rm and} \; \lambda_1 + \cdots+ \lambda_r \geq
     \mu_1 + \cdots + \mu_r \; {\rm for} \;{\rm all}\; r\geq1.
\ee
A partition is identified with the Young diagram.
The conjugate partition of $\lambda$,
whose diagram is obtained by interchanging rows and columns,
is denoted by $\lambda'$.
$x^{\lambda}$ stands for the monomial
$x^{\lambda}=x_1^{\lambda_1}x_2^{\lambda_2}\cdots$.
We give some bases of $\Lambda_{{\bf Q}(q,t)}$ :\\
(\romannumeral1)
$\displaystyle
  m_{\lambda}=
  \sum_{\mbox{\scriptsize $\alpha$ : distinct} \atop
        \mbox{\scriptsize permutation of $\lambda$}}
  \!\!\!\!x^{\alpha},\quad
(\mbox{monomial symmetric function})$. \\
(\romannumeral2)
$\displaystyle
  p_{\lambda}=p_{\lambda_1}p_{\lambda_2}\cdots,\quad
  (\mbox{the $r$-th power sum } p_r=\sum_i x_i^r)$. \\
(\romannumeral3)
$\displaystyle
  e_{\lambda}=e_{\lambda_1}e_{\lambda_2}\cdots,\quad
  (\mbox{the $r$-th elementary symmetric function }
   e_r=\!\sum_{i_1<\cdots<i_r} x_{i_1}\cdots x_{i_r})$.

We endow $\Lambda_{{\bf Q}(q,t)}$ with the following
scalar product:
\begin{equation}
\langle p_{\lambda},p_{\mu}\rangle_{q,t} =
 \delta_{\lambda,\mu}z_{\lambda}(q,t),
\end{equation}
where
\begin{eqnarray}
z_{\lambda}(q,t)=\prod_{r\geq1} r^{m_r} m_r ! \cdot
  \prod_{i=1}^{\ell(\lambda)}
\frac{1-q^{\lambda_i}}{1-t^{\lambda_i}},
\quad \lambda=(1^{m_1}2^{m_2}\cdots).
\end{eqnarray}
with $m_r\equiv \#\{i\,|\,\lambda_i=r\}$.

The Macdonald symmetric function is
characterized by the following existence theorem:
\begin{theorem} {\rm \cite{rMac}}
For each partition $\lambda$ there is a unique symmetric function
$P_{\lambda}=P_{\lambda}(x;q,t)\in \Lambda_{{\bf Q}(q,t)}$ such that
\ba
  {\rm (A)}
  &&
  P_{\lambda}=m_{\lambda} + \sum_{\mu<\lambda} u_{\lambda\mu}m_{\mu},
  \quad
  u_{\lambda\mu}\in {\bf Q}(q,t), \\
  {\rm (B)}
  &&
  \langle P_{\lambda},P_{\mu}\rangle_{q,t} =0\qquad { \rm if } \;
  \lambda\neq \mu.
\ea
\end{theorem}

Even though this definition is concise, it is more
useful to define the Macdonald symmetric function by
introducing an operator which can be regarded as a natural
deformation of the Calogero-Sutherland Hamiltonian.
The following operator is called as the Macdonald operator:
\begin{equation}
D_1 = \sum_{i=1}^n \left( \prod_{j\neq i}
         \frac{tx_i-x_j}{x_i-x_j}\right)T_{q,x_i},
\end{equation}
where $T_{q,x_i}$ is the $q$-shift operator defined by $
(T_{q,x_i}f)(x_1,\cdots,x_n)=f(x_1,\cdots,qx_i,\cdots,x_n)$.
The other way to define the Macdonald symmetric functions is the
following:
\begin{theorem} {\rm \cite{rMac}}
For each partition $\lambda$
$($of length $\leq n$$)$, there is a unique symmetric function
$P_{\lambda}(x;q,t)\in\Lambda_{n,{\bf Q}(q,t)}$
satisfying the two conditions:
\ba
  {\rm (A)}
  &&
  P_{\lambda}=m_{\lambda} + \sum_{\mu<\lambda} u_{\lambda\mu}m_{\mu},
  \quad
  u_{\lambda\mu}\in {\bf Q}(q,t), \\
  {\rm (C)}
  &&
  D_1 P_{\lambda}=\sum_{i=1}^n t^{n-i}q^{\lambda_i}\cdot P_{\lambda}.
\ea
\end{theorem}

It was shown by Macdonald that $D_1$ can be included in a family
of mutually commuting operators $\{D_r|r=1,\cdots,n\}$.
The operator $D_r$ is defined by
\begin{equation}
D_r=\sum_{I\subset \{1,2,\cdots,n\}\atop |I|=r}
t^{r(r-1)/2}\prod_{i\in I \atop j \not\in I}
      \frac{tx_i-x_j}{x_i-x_j}\cdot
  \prod_{i\in I} T_{q,x_i}.
\end{equation}
\begin{theorem} {\rm \cite{rMac}}
The operators $D_r$ $($$r=1,\cdots,n$$)$ commute with each other and
the Macdonald symmetric function
$P_{\lambda}$ is the simultaneous eigenvector of these operators
with the eigenvalues given by the coefficient of $X^{n-r}$
in $\prod_{i=1}^n (X+t^{n-i}q^{\lambda_i})$.
\end{theorem}
One can notice that Macdonald's $D_r$ and Ruisenaars' $\hat{S}_k$
have similar structure.

Here we list some particular cases of the Macdonald symmetric
function $P_{\lambda}(q,t)$.\\
(\romannumeral1)
 When $t=q$, $P_{\lambda}(q,q)$ is the Schur function $s_{\lambda}$.\\
(\romannumeral2)
 When $q=0$, $P_{\lambda}(0,t)$ is the Hall-Littlewood function
 $P_{\lambda}(t)$.\\
(\romannumeral3)
$\displaystyle
 \lim_{q \rightarrow 1} P_{\lambda}(q,q^{\beta})$
 is the Jack symmetric function $J_{\lambda}(\beta)$.\\
(\romannumeral4)
 When $t=1$, $P_{\lambda}(q,1)=m_{\lambda}$.\\
(\romannumeral5)
 When $q=1$, $P_{\lambda}(1,t)=e_{\lambda'}$.\\
(\romannumeral6)
 $P_{\lambda}(q^{-1},t^{-1})=P_{\lambda}(q,t)$.

Let $Q_{\lambda}$'s be the dual basis of $P_{\lambda}$'s, that is
\begin{equation}
\langle Q_{\lambda}(q,t),P_{\mu}(q,t)\rangle_{q,t}=\delta_{\lambda,\mu}.
\end{equation}
The following proposition is easily proved:
\begin{prop}
We have the following equation:
\begin{equation}
\sum_{\lambda} P_{\lambda}(x;q,t) Q_{\lambda}(y;q,t)
= \Pi(x,y;q,t),
\end{equation}
where
\begin{equation}
\Pi(x,y)=\Pi(x,y;q,t)=\prod_{i,j}
      \frac{(tx_i y_j;q)_{\infty}}{(x_i y_j;q)_{\infty}}.
\end{equation}
\end{prop}
Here we have used the following notation:
\begin{equation}
(a;q)_{\infty}=\prod_{k=0}^{\infty}(1-aq^k)\qquad {\rm for}\;\;
 a\in{\bf C}.
\end{equation}
We remark that $\Pi(x,y)=\Pi(y,x)$ is a Taylor series in $x_i$
and $y_j$ in the region $|x_iy_j|<1$.

If we write
\begin{equation}
Q_{\lambda}(x;q,t)= b_{\lambda}(q,t)P_{\lambda}(x;q,t),
\end{equation}
then we have the explicit formula for $b_{\lambda}(q,t)$.
To state it we need the following notation;
For any $s=(i,j)\in\lambda$ ($i$-th row, $j$-th column in Young
diagram $\lambda$),
let us define arm-length $a(s)$, leg-length $\ell(s)$,
arm-colength $a'(s)$ and leg-colength $\ell'(s)$ as follows:
\begin{equation}
\left\{
\begin{array}{cc}
a(s)=\lambda_i-j, & a'(s)=j-1,\\
\ell(s)=\lambda'_j-i, & \ell'(s)=i-1.
\end{array} \right.
\end{equation}
\begin{theorem} {\rm \cite{rMac}}
The explicit formula for the coefficient $b_{\lambda}(q,t)$ is
\begin{equation}
   b_{\lambda}(q,t)=
   \frac{1}{\langle P_{\lambda},P_{\lambda}\rangle_{q,t}}=
   \prod_{s\in\lambda} \frac{1-q^{a(s)}t^{\ell(s)+1}}
                            {1-q^{a(s)+1}t^{\ell(s)}}.
\end{equation}
\end{theorem}

\noindent
{\large\bf 2.2}~{\sl The dual transformation}\\
Let us define an automorphism
\begin{equation}
\omega_{q,t}:\Lambda_{{\bf Q}(q,t)}\rightarrow
          \Lambda_{{\bf Q}(q,t)}
\end{equation}
by fixing the action on $p_r$ as
\begin{equation}
\omega_{q,t}(p_r)=(-1)^{r-1}\frac{1-q^r}{1-t^r}p_r
\end{equation}
and extending it naturally.
We have the following theorem which describes the duality
transformation of the Macdonald symmetric functions.
\begin{theorem} {\rm \cite{rMac}}
For any partition $\lambda$, we have
\begin{equation}
\omega_{q,t}P_{\lambda}(q,t)=Q_{\lambda'}(t,q) \qquad
\Bigl(\mbox{\rm or equivalently }
\omega_{q,t}Q_{\lambda}(q,t)=P_{\lambda'}(t,q) \Bigr).
\end{equation}
\end{theorem}

It is easy to show 
\begin{equation}
\omega_{q,t}^y \Pi(x,y;q,t)=\prod_{i,j}(1+x_i y_j)\equiv
\tilde{\Pi}(x,y),
\end{equation}
where $\omega_{q,t}^y$ acts on the variable $y$.
Hence we have
\begin{equation}
\sum_{\lambda} P_{\lambda}(x;q,t) P_{\lambda'}(y;t,q)
= \tilde{\Pi}(x,y)
\end{equation}

\noindent
{\large\bf 2.3}~{\sl The scalar product
$\langle\;,\;\rangle'_{n;q,t}$}\\
Next we consider the properties of
another scalar product $\langle~,~\rangle'_{n;q,t}$
that will be defined below.
We shall work with a finite number of indeterminates
$x=(x_1,\cdots,x_n)$.
We set the parameters $q$ and $t$ as $0<q<1$ and $0<t<1$.
Define
\begin{equation}
\Delta(x)=\Delta(x;q,t)=\prod_{i,j=1 \atop i\neq j}^n
  \frac{(x_i/x_j;q)_{\infty}}{(t x_i/x_j;q)_{\infty}}.
\end{equation}
In the region $t<|x_i/x_j|<t^{-1}$ ($i\neq j$),
$\Delta$ is a Laurent series in $x_i$'s.
For $f,g \in \Lambda_{n,{\bf Q}(q,t)}$, we define\footnote{
For $f(x)=f(x_1,x_2,\cdots)$,
we define $f(\overline{x})=f(1/x_1,1/x_2,\cdots)$.}
\ba
  \langle f,g \rangle'_{n;q,t}
  &\!\!=\!\!&
  \frac{1}{n!} \oint \prod_{j=1}^n \frac{dx_j}{2\pi i x_j}\cdot
  f(\overline{x})g(x)\Delta(x;q,t) \\
  &\!\!=\!\!&
  \frac{1}{n!} \Bigl(\mbox{constant term in }
  f(\overline{x})g(x)\Delta(x)\Bigr).
  \nonumber
\ea
The following proposition is the most fundamental one.
\begin{prop} {\rm \cite{rMac}}
The operator $D_1$ is self-adjoint
with respect to this scalar product, namely,
\begin{equation}
\langle D_1 f,g \rangle'_{n;q,t}=\langle f,D_1 g \rangle'_{n;q,t} \:,
\end{equation}
for all  $f,g \in \Lambda_{n,{\bf Q}(q,t)}$.
\end{prop}
{}From this proposition we have
\begin{prop} {\rm \cite{rMac}}
\begin{equation}
\langle P_{\lambda}(q,t),P_{\mu}(q,t)\rangle'_{n;q,t}=0\qquad
{\rm if}\;\lambda\neq\mu.
\end{equation}
\end{prop}
Furthermore we have the following conjecture:
\begin{conj}{\bf (Macdonald's constant term conjecture)}
{\rm \cite{rMac}}
\begin{eqnarray}
\langle P_{\lambda}(q,t),P_{\lambda}(q,t)\rangle'_{n;q,t}&\!\!=\!\!&
\prod_{1\leq i<j\leq n}
 \frac{(q^{\lambda_i-\lambda_j}t^{j-i};q)_{\infty}
       (q^{\lambda_i-\lambda_j+1}t^{j-i};q)_{\infty}}
      {(q^{\lambda_i-\lambda_j}t^{j-i+1};q)_{\infty}
       (q^{\lambda_i-\lambda_j+1}t^{j-i-1};q)_{\infty}}\\
&\!\!=\!\!& \prod_{i=1}^n \frac{\Gamma_q(i\beta)}
                    {\Gamma_q(\beta)\Gamma_q((i-1)\beta+1)}\cdot
     \prod_{s\in\lambda}\frac{1-q^{a'(s)}t^{n-\ell'(s)}}
                             {1-q^{a'(s)+1}t^{n-\ell'(s)-1}}\cdot
     b^{-1}_{\lambda}(q,t),\nonumber
\end{eqnarray}
where $\Gamma_q(x)$ is the $q$-gamma function defined by
$$
\Gamma_q(x)=\frac{(q;q)_{\infty}}{(q^x;q)_{\infty}}(1-q)^{1-x}.
$$
\end{conj}

\section{Integral Representation of the Macdonald Symmetric Functions}
\setcounter{equation}{0}

In this section, we construct integral representation
formulas for the Macdonald symmetric functions.
We adopt the same idea as that of
the case of Jack symmetric functions\cite{rAMOS1,rMY1,rMY2,rAMOS2}.

\noindent
{\large\bf 3.1}~{\sl Maps ${\cal G}_s$, ${\cal N}_{n,m}$ and
an integral formula for $P_{\lambda}(x;q,t)$}\\
Let us define the map ${\cal G}_s$ and ${\cal N}_{n,m}$ as follows:
\begin{eqnarray}
  {\cal G}_s
  &:&
  \Lambda_{r,{\bf Q}(q,t)}\rightarrow \Lambda_{r,{\bf Q}(q,t)}\\
  &&
  f(x)\mapsto({\cal G}_sf)(x)=\prod_{i=1}^r (x_i)^s\cdot f(x),
  \n
  {\cal N}_{n,m}
  &:&
  \Lambda_{m,{\bf Q}(q,t)}\rightarrow \Lambda_{n,{\bf Q}(q,t)} \\
  &&
  f(x)\mapsto({\cal N}_{n,m}f)(x)=
  \oint\prod_{j=1}^m\frac{dy_j}{2\pi iy_j}\cdot
  \Pi(x,\overline{y};q,t)\Delta(y;q,t)f(y).
  \nonumber
\end{eqnarray}
Here $r,m<\infty$ and $n$ can be equal to $\infty$.

\begin{prop}The actions of ${\cal G}_s$ and
${\cal N}_{n,m}$ on the Macdonald symmetric function
$P_{\lambda}$ are as follows:
\begin{eqnarray}
  &{\rm (\romannumeral1)}&
  P_{(s^r)+\lambda}(x_1,\cdots,x_r;q,t)=
  {\cal G}_sP_{\lambda}(x_1,\cdots,x_r;q,t), \\
  &{\rm (\romannumeral2)}&
  P_{\lambda}(x_1,\cdots,x_n;q,t)=
  \frac{\langle P_{\lambda},P_{\lambda}\rangle_{q,t}}
       {m! \langle P_{\lambda},P_{\lambda}\rangle'_{m;q,t}}
  {\cal N}_{n,m}P_{\lambda}(x_1,\cdots,x_m;q,t).
\end{eqnarray}
\end{prop}

{\it Proof} \quad As for (\romannumeral1),
we can easily check the conditions (A) and (C) in theorem 2.2.
The statement (\romannumeral2) can be proved as follows:
\begin{eqnarray*}
  &&
  \frac{\langle P_{\lambda},P_{\lambda}\rangle_{q,t}}
       {m! \langle P_{\lambda},P_{\lambda}\rangle'_{m;q,t}}
  \oint \prod_{j=1}^m\frac{dy_j}{2\pi i y_j}\cdot
  \Pi(x,\overline{y};q,t)\Delta(y;q,t)P_{\lambda}(y;q,t) \\
  &=\!\!&
  \frac{\langle P_{\lambda},P_{\lambda}\rangle_{q,t}}
       {m! \langle P_{\lambda},P_{\lambda}\rangle'_{m;q,t}}
  \oint \prod_{j=1}^m\frac{dy_j}{2\pi i y_j}\cdot
  \sum_{\mu}Q_{\mu}(x;q,t)P_{\mu}(\overline{y};q,t)
  \Delta(y;q,t)P_{\lambda}(y;q,t) \\
  &=\!\!&
  \langle P_{\lambda},P_{\lambda}\rangle_{q,t}
  Q_{\lambda}(x;q,t)\\
  &=\!\!&
  P_{\lambda}(x;q,t).
\qquad\qquad\qquad\qquad\qquad\qquad\qquad\qquad
\qquad\qquad\qquad\qquad
  Q.E.D.
\end{eqnarray*}

Any Young diagram $\lambda$ can be uniquely decomposed into
rectangles:
\generalYoung
\noindent
where $r_N>\cdots>r_2>r_1$.
Therefore the partition $\lambda$ is parametrized as follows:
\begin{equation}
  \lambda=(s_N^{r_N})+\cdots+(s_2^{r_2})+(s_1^{r_1}),
  \label{deflam}
\end{equation}
where $(s^r)=(\overbrace{s,s,\cdots,s}^r)$.
For the partition $\lambda$, we assign a set of partitions
$\lambda^{(a)}$ ($a=1,\cdots,N$)
as follows:
\be
  \lambda^{(a)}=(s_a^{r_a})+\cdots+(s_2^{r_2})+(s_1^{r_1}).
\ee

Here we state our main theorem:
\begin{theorem}
Let $\lambda$ be the partition given by {\rm (\ref{deflam})}.
We have the following multiple integral representation of
the Macdonald symmetric function
$P_{\lambda}(x;q,t)\in\Lambda_{{\bf Q}(q,t)}$,
\begin{eqnarray}
  P_{\lambda}(x;q,t)
  &\!\!=\!\!&
  C^+_{\lambda}
  {\cal N}_{r_{N+1},r_N}{\cal G}_{s_N}{\cal N}_{r_N,r_{N-1}}\cdots
  {\cal G}_{s_2}{\cal N}_{r_2,r_1}{\cal G}_{s_1}\cdot 1
  \n
  &\!\!=\!\!&
  C^+_{\lambda}
  \oint\prod_{a=1}^N\prod_{j=1}^{r_a}
  \frac{dx^a_j}{2 \pi ix^a_j}\Bigl(x^a_j\Bigr)^{s_a}\cdot
  \prod_{a=1}^N\Pi(x^{a+1},\overline{x^a};q,t)\Delta(x^a;q,t),
  \label{intmac}
\end{eqnarray}
where $x_i=x^{N+1}_i$, $r_{N+1}=\infty$ and
\begin{equation}
  C^+_{\lambda}=C^+_{\lambda}(q,t)=
  \prod_{a=1}^N
  \frac{\langle P_{\lambda^{(a)}},P_{\lambda^{(a)}}\rangle_{q,t}}
    {r_a!\langle P_{\lambda^{(a)}},P_{\lambda^{(a)}}\rangle'_{r_a;q,t}}.
\end{equation}
\end{theorem}
{\it Proof}\qquad Use proposition 3.1 iteratively.
${\cal G}_s$ adds a rectangle and ${\cal N}_{n,m}$ increases the
number of variables. \qquad{\it Q.E.D.}

\noindent
{\large\bf 3.2}~{\sl Another integral formula for
$P_{\lambda'}(x;t,q)$}\\
Next we consider
another integral representation of the Macdonald symmetric function
$P_{\lambda'}(x;t,q)$ that is obtained from
$Q_{\lambda}(x;q,t)$ by applying the automorphism $\omega_{q,t}$.
Let us introduce one more map defined by
\begin{eqnarray}
  \tilde{\cal N}_{n,m}
  &:&
  \Lambda_{m,{\bf Q}(q,t)}\rightarrow\Lambda_{n,{\bf Q}(q,t)}\\
  &&
  f(x) \mapsto
  (\tilde{\cal N}_{n,m}f)(x)=
  \oint\prod_{j=1}^m\frac{dy_j}{2\pi i y_j}\cdot
  \tilde{\Pi}(x,\overline{y})\Delta(y;q,t)f(y).
  \nonumber
\end{eqnarray}
We can prove the following proposition 3.3 and theorem 3.4
in the same way as proposition 3.1 and theorem 3.2, respectively.
\begin{prop}
The following equality holds:
\begin{equation}
  P_{\lambda'}(x_1,\cdots,x_n;t,q)=
  \frac{1}{m!\langle P_{\lambda},P_{\lambda}\rangle'_{m;q,t}}
  \tilde{\cal N}_{n,m}P_{\lambda}(x_1,\cdots,x_m;q,t).
\end{equation}
\end{prop}

\begin{theorem}
Let $\lambda$ be the partition given by {\rm (\ref{deflam})}.
We have the following multiple integral representation of
the Macdonald symmetric function
$P_{\lambda'}(x;t,q)\in\Lambda_{{\bf Q}(q,t)}$,
\begin{eqnarray}
  P_{\lambda'}(x;t,q)
  &\!\!=\!\!&
  C^-_{\lambda}
  \tilde{\cal N}_{r_{N+1},r_N}{\cal G}_{s_N}{\cal N}_{r_N,r_{N-1}}\cdots
  {\cal G}_{s_2}{\cal N}_{r_2,r_1}{\cal G}_{s_1}\cdot 1
  \n
  &\!\!=\!\!&
  C^-_{\lambda}
  \oint\prod_{a=1}^N\prod_{j=1}^{r_a}
  \frac{dx^a_j}{2 \pi i x^a_j}\Bigl(x^a_j\Bigr)^{s_a}\\
  &&
  \qquad\times\:
  \tilde{\Pi}(x^{N+1},\overline{x^N};q,t)
  \prod_{a=1}^{N-1}\Pi(x^{a+1},\overline{x^a};q,t)\cdot
  \prod_{a=1}^N\Delta(x^a;q,t),
  \nonumber
\end{eqnarray}
where $x_i=x^{N+1}_i$, $r_{N+1}=\infty$ and
\begin{equation}
  C^-_{\lambda}=C^-_{\lambda}(q,t)=
  \frac{C^+_{\lambda}(q,t)}
       {\langle P_{\lambda},P_{\lambda}\rangle_{q,t}}.
\end{equation}
\end{theorem}

\noindent
{\large\bf 3.3}~{\sl An integral formula for the skew Macdonald functions}\\
Now let us proceed to discuss how to obtain an integral
representation of the skew Macdonald functions. To this end,
let us start with introducing a boson Fock space $\cal{F}$
which is isomorphic to the ring of the
symmetric functions $\Lambda_{{\bf Q}(q,t)}$ \cite{rDJKM,rAMOS1}.
Define the commutation relations of the bosonic oscillators
$a_n$ ($n\in{\bf Z}_{\neq 0}$) as follows:
\begin{equation}
[a_n,a_m]= n \frac{1-q^{|n|}}{1-t^{|n|}}\delta_{n+m,0}.
\end{equation}
Let $|0\rangle$ be the vacuum vector such that
$a_n|0\rangle=0$ for $n<0$ and ${\cal F}$ be the
Fock space defined by ${\cal F}=
{\bf Q}(q,t)[a_{-1},a_{-2},\cdots]|0\rangle$.
Let $\langle 0|$ be the dual of $|0\rangle$ i.e.,
$\langle 0|0\rangle=1$. Define
${\cal F}^*=\langle 0|{\bf Q}(q,t)[a_{1},a_{2},\cdots]$.

We can construct an isomorphism $\iota$ between
${\cal F}$ and $\Lambda_{{\bf Q}(q,t)}$ as follows:
\begin{eqnarray}
\iota &:& {\cal F}\rightarrow \Lambda_{{\bf Q}(q,t)}\\
      && |f\rangle \mapsto f(x)=\langle 0|C(x)|f\rangle,
\end{eqnarray}
and an isomorphism $\iota^*$ between
${\cal F}^*$ and $\Lambda_{{\bf Q}(q,t)}$ by
\begin{eqnarray}
\iota^* &:& {\cal F}^*\rightarrow \Lambda_{{\bf Q}(q,t)}\\
      && \langle f|  \mapsto f(x)=\langle f| C^*(x) |0\rangle,
\end{eqnarray}
where
\begin{eqnarray}
C(x) &\!\!=\!\!&   \exp\left( \sum_{n=1}^\infty \frac{1-t^{n}}{1-q^n}
          \frac{a_{n} }{n} p_n\right),\\
C^*(x) &\!\!=\!\!&   \exp\left( \sum_{n=1}^\infty \frac{1-t^{n}}{1-q^n}
          \frac{a_{-n} }{n} p_n\right).
\end{eqnarray}
We recall that $p_n$ is the power sum $p_n=x_1^n+x_2^n+\cdots$.

We will use the following notation. For any symmetric function
$f\in\Lambda_{{\bf Q}(q,t)}$,
we assign an operator $\hat{f}\in {\bf Q}(q,t)[a_{-1},a_{-2},\cdots]$
and a vector $|f\rangle\in{\cal F}$ such that
$\iota(\hat{f}|0\rangle)=\iota(|f\rangle)=f(x)$.
In the same way,
we assign an operator $\hat{f}^*\in {\bf Q}(q,t)[a_{1},a_{2},\cdots]$
and a vector $\langle f|\in{\cal F}^*$ such that
$\iota^*(\langle 0|\hat{f}^*)=\iota^*(\langle f |)=f(x)$. For example,
$\iota(\hat{P}_{\lambda}(q,t)|0\rangle)=\iota(|P_{\lambda}(q,t)\rangle)=
P_{\lambda}(x;q,t)$ and
$\iota^*(\langle 0|\hat{P}^*_{\lambda}(q,t))=
\iota^*(\langle P_{\lambda}(q,t)|)=P_{\lambda}(x;q,t)$.
For a product $f(x)g(x)$, the corresponding state is
$\hat{f}\hat{g}\ket{0}=\hat{f}\ket{g}=\ket{f\cdot g}$.
We have the following proposition:
\begin{prop}~\\
{\rm (\romannumeral1)}
Let $\bra{f}\in{\cal F}^*$ and $\ket{g}\in{\cal F}$.
We have
\begin{equation}
\langle f |g\rangle=\langle f(x),g(x)\rangle_{q,t}.
\end{equation}
{\rm (\romannumeral2)}
Let $\langle f | \in{\cal F}^*$ and $|g\cdot h\rangle\in{\cal F}$.
We have
\begin{equation}
\langle f|g\cdot h\rangle=
\langle \langle f|C^*(x)|g\rangle,\langle 0|C(x)|h\rangle \rangle_{q,t}.
\end{equation}
\end{prop}

{\it Proof.}\qquad
We defined the commutation relations of $a_n$ so that
(\romannumeral1) holds.\\
(\romannumeral2) is proved as follows:
\begin{eqnarray*}
  \langle f|g\cdot h\rangle
  &\!\!=\!\!&
  \langle f|\hat{g}|h\rangle
  =\langle i\ket{h}\\
  &\!\!=\!\!&
  \langle \langle i|C^*(x)|0\rangle,
  \langle 0|C(x)|h\rangle \rangle_{q,t}=
  \langle \langle f|\hat{g}C^*(x)|0\rangle,
  \langle 0|C(x)|h\rangle \rangle_{q,t}\\
  &\!\!=\!\!&
  \langle \langle f|C^*(x)|g\rangle,
  \langle 0|C(x)|h\rangle \rangle_{q,t},
\end{eqnarray*}
where we have set $\bra{i}=\bra{f}\hat{g}\in{\cal F}^*$
(which may be $0$), and used (\romannumeral1) and
$\hat{g}C^*(x)=C^*(x)\hat{g}$.
\qquad{\it Q.E.D.}\\
We remark that in this boson language, for example,
proposition 2.4 is a consequence of the completeness condition
$\sum_{\lambda}\ket{P_{\lambda}}\bra{Q_{\lambda}}={\bf 1}$.

By theorems 3.2 and 3.4,
we have the following bosonization formulas for
the Macdonald symmetric function $P_{\lambda}(x;q,t)$:
\begin{prop}
Let $\lambda$ be the partition given by {\rm (\ref{deflam})}.
We have
\begin{eqnarray}
  \hat{P}_{\lambda}(q,t)
  &\!\!=\!\!&
  \oint \prod_{j=1}^{r_N} \frac{dx_j}{2\pi i x_j}\cdot
  F^+_{\lambda}(x;q,t)
  \prod_{j=1}^{r_N}
  \exp\left( \sum_{n=1}^\infty \frac{1-t^n}{1-q^n}
  \frac{a_{-n} }{n} (x_j)^{-n}\right) \\
  &\!\!=\!\!&
  \oint \prod_{j=1}^{r_N} \frac{dx_j}{2\pi i x_j}\cdot
  F^-_{\lambda}(x;q,t)
  \prod_{j=1}^{r_N}
  \exp\left( -\sum_{n=1}^\infty
  \frac{a_{-n} }{n} (-x_j)^{-n}\right), \\
  \hat{P}_{\lambda}^*(q,t)
  &\!\!=\!\!&
  \oint \prod_{j=1}^{r_N} \frac{dx_j}{2\pi i x_j}\cdot
  F^+_{\lambda}(x;q,t)
  \prod_{j=1}^{r_N}
  \exp\left( \sum_{n=1}^\infty \frac{1-t^n}{1-q^n}
  \frac{a_{n} }{n} (x_j)^{-n}\right) \\
  &\!\!=\!\!&
  \oint \prod_{j=1}^{r_N} \frac{dx_j}{2\pi i x_j}\cdot
  F^-_{\lambda}(x;q,t)
  \prod_{j=1}^{r_N}
  \exp\left( -\sum_{n=1}^\infty
  \frac{a_{n} }{n} (-x_j)^{-n}\right).
\end{eqnarray}
Namely,
$\langle 0| C(x)\hat{P}_{\lambda}(q,t)|0\rangle=
\langle 0| \hat{P}^*_{\lambda}(q,t)C^*(x)|0\rangle=
P_{\lambda}(x;q,t)$.
Here $F^{\pm}_{\lambda}$ is defined by
\ba
  F^{+}_{\lambda}(x^N;q,t)
  &\!\!=\!\!&
  C^+_{\lambda}(q,t)
  \Delta(x^N;q,t)\prod_{j=1}^{r_N}\Bigl(x^N_j\Bigr)^{s_N} \n
  &&\times
  \oint\prod_{a=1}^{N-1}\prod_{j=1}^{r_a}
  \frac{dx^a_j}{2\pi ix^a_j}\Bigl(x^a_j\Bigr)^{s_a}\cdot
  \prod_{a=1}^{N-1}\Pi(x^{a+1},\overline{x^a};q,t)\Delta(x^a;q,t),\\
  F^{-}_{\lambda'}(x^N;q,t)
  &\!\!=\!\!&
  C^-_{\lambda}(t,q)
  \Delta(x^N;t,q)\prod_{j=1}^{r_N}\Bigl(x^N_j\Bigr)^{s_N} \n
  &&\times
  \oint\prod_{a=1}^{N-1}\prod_{j=1}^{r_a}
  \frac{dx^a_j}{2\pi ix^a_j}\Bigl(x^a_j\Bigr)^{s_a}\cdot
  \prod_{a=1}^{N-1}\Pi(x^{a+1},\overline{x^a};t,q)\Delta(x^a;t,q).
\ea
\end{prop}


Let $\mu$ and $\nu$ be two partitions. We define the structure constants
$f^{\lambda}_{\mu\nu}$ of the
ring $\Lambda_{{\bf Q}(q,t)}$ by
\begin{equation}
  P_{\mu}(x;q,t)P_{\nu}(x;q,t)=
  \sum_{\lambda}f^{\lambda}_{\mu\nu}(q,t)P_{\lambda}(x;q,t),
\end{equation}
or equivalently,
\begin{equation}
f^{\lambda}_{\mu\nu}=f^{\lambda}_{\mu\nu}(q,t)=
\langle Q_{\lambda},P_{\mu}P_{\nu}\rangle_{q,t}\in {\bf Q}(q,t).
\end{equation}
The skew $Q$-function is defined by
\begin{equation}
  Q_{\lambda/\mu}(x;q,t)=
  \sum_{\nu}f^{\lambda}_{\mu\nu}(q,t)Q_{\nu}(x;q,t).
\end{equation}
This is equivalent to the following condition:
\begin{equation}
\langle Q_{\lambda/\mu},P_{\nu} \rangle_{q,t}=
\langle Q_{\lambda},P_{\mu}P_{\nu} \rangle_{q,t}.
\end{equation}
The skew $P$-function is given by
$P_{\lambda/\mu}=b_{\lambda}^{-1}b_{\mu}Q_{\lambda/\mu}$.

Now we are ready to state the boson representation of
the skew $Q$-function.

\begin{theorem}We have the following boson realization
of the skew $Q$-function.
\begin{equation}
Q_{\lambda/\mu}(x;q,t)=
\langle Q_{\lambda}|C^*(x)|P_{\mu}\rangle
=\langle P_{\mu}|C(x)|Q_{\lambda}\rangle .
\end{equation}
\end{theorem}

{\it Proof.}\qquad From propositions 3.5 and 3.6 we have the following:
\begin{eqnarray*}
\langle \langle Q_{\lambda}|C^*(x)|P_{\mu}\rangle  ,P_{\nu}\rangle_{q,t}
&\!\!=\!\!&\langle \langle Q_{\lambda}|C^*(x)|P_{\mu}\rangle ,
   \langle 0|C(x)|P_{\nu}\rangle \rangle_{q,t}\\
&\!\!=\!\!&\langle  Q_{\lambda},P_{\mu}P_{\nu}\rangle_{q,t}.
\end{eqnarray*}
This proves the first equality. The second equality can be proved
in the same way.\qquad{\it Q.E.D.}

As a corollary of this theorem and proposition 3.6,
we obtain integral representation
formulas for the skew $Q$-function.
\begin{cor}
Let $\lambda$ be the partition given by {\rm (\ref{deflam})}
and $\mu$ be another partition
$\mu=(\sigma_M^{\rho_M})+\cdots+(\sigma_1^{\rho_1})$.
We have the integral representation formulas for
$ Q_{\lambda/\mu}(x;q,t)\in \Lambda_{{\bf Q}(q,t)}$ as follows:
\begin{eqnarray}
  &&
  b(q,t)_{\lambda}^{-1}Q_{\lambda/\mu}(x;q,t)\n
  &=\!\!&
  \oint\prod_{j=1}^{r_N}\frac{dz_j}{2\pi i z_j} \cdot
       \prod_{j=1}^{\rho_M}\frac{dw_j}{2\pi i w_j}\cdot
  F^+_{\lambda}(z;q,t)F^+_{\mu}(w;q,t)
  \Pi(\overline{z},\overline{w};q,t)
  \times\left\{
  \begin{array}{l}
    \Pi(x,\overline{z};q,t) \\
    \Pi(x,\overline{w};q,t)
  \end{array}\right. \\
  &=\!\!&
  \oint\prod_{j=1}^{r_N}\frac{dz_j}{2\pi i z_j} \cdot
       \prod_{j=1}^{\rho_M}\frac{dw_j}{2\pi i w_j}\cdot
  F^-_{\lambda}(z;q,t)F^+_{\mu}(w;q,t)
  \tilde{\Pi}(\overline{z},\overline{w})
  \times\left\{
  \begin{array}{l}
    \tilde{\Pi}(x,\overline{z}) \\
    \Pi(x,\overline{w};q,t)
  \end{array}\right. \\
  &=\!\!&
  \oint\prod_{j=1}^{r_N}\frac{dz_j}{2\pi i z_j} \cdot
       \prod_{j=1}^{\rho_M}\frac{dw_j}{2\pi i w_j}\cdot
  F^+_{\lambda}(z;q,t)F^-_{\mu}(w;q,t)
  \tilde{\Pi}(\overline{z},\overline{w})
  \times\left\{
  \begin{array}{l}
    \Pi(x,\overline{z};q,t) \\
    \tilde{\Pi}(x,\overline{w})
  \end{array}\right. \\
  &=\!\!&
  \oint\prod_{j=1}^{r_N}\frac{dz_j}{2\pi i z_j} \cdot
       \prod_{j=1}^{\rho_M}\frac{dw_j}{2\pi i w_j}\cdot
  F^-_{\lambda}(z;q,t)F^-_{\mu}(w;q,t)
  \Pi(\overline{z},\overline{w};t,q)
  \times\left\{
  \begin{array}{l}
    \tilde{\Pi}(x,\overline{z}) \\
    \tilde{\Pi}(x,\overline{w}).
  \end{array}\right.
\end{eqnarray}
\end{cor}


{\it Remark.}\quad
More generally, one can directly prove that
the skew Macdonald polynomial can be written
in the integral transformation ${\cal N}_{n,m}$ of eq.\ (3.2) or
in the power-sums $p_n$ as follows:
\ba
  Q_{\lambda/\mu}(x;q,t)
  &\!\!=\!\!&
  {\langle Q_\lambda,Q_\lambda\rangle_{q,t} \over
 m!\langle Q_\lambda,Q_\lambda\rangle'_{m;q,t} }
 \left({\cal N}_{n,m}   Q_\lambda\,\overline{P}_\mu \right)(x;q,t), \n
  Q_{\lambda/\mu}(p\,;q,t)
  &\!\!=\!\!&
  P_\mu(\overline{p}\,;q,t)\,\, Q_\lambda(p;q,t) \cdot 1,
\ea
for all $m\geq\ell(\lambda)$.
Here $\overline{P_\mu(x)}\equiv P_\mu\left({1\over x}\right)$ and
$\overline p_n\equiv n{1-q^n\over 1-t^n}{\partial\over\partial p_n}$.

\noindent
{\large\bf 3.4}~{\sl The Kostka matrix}\\
As another application of the bosonization constructed in the last
subsection, we will give integral representations of the Kostka
matrix $K_{\lambda\mu}(q,t)$.
Let
\begin{eqnarray}
h_{\lambda}(q,t)&=& \prod_{s\in \lambda} (1-q^{a(s)}t^{\ell(s)+1}),\\
h_{\lambda}'(q,t)&=& \prod_{s\in \lambda}
(1-q^{a(s)+1}t^{\ell(s)})=h_{\lambda}(t,q).
\end{eqnarray}
So, we have
\begin{equation}
b_{\lambda}(q,t)=h_{\lambda}(q,t)/h_{\lambda}'(q,t).
\end{equation}
Let us define
\be
  M_{\lambda}(x;q,t)=
  h_{\lambda}(q,t)P_{\lambda}(x;q,t)=
  h_{\lambda}'(q,t)Q_{\lambda}(x;q,t).
\ee
The $q$-analogue of the Kostka-Foulks polynomial $K_{\lambda\mu}(q,t)$
introduced by Macdonald is defined by
\begin{equation}
M_{\mu}(x;q,t)=\sum_{\lambda}K_{\lambda\mu}(q,t)S_{\lambda}(x;t),
\end{equation}
where $S_{\lambda}(x;t)$ is the dual base of the Schur function
$s_{\lambda}(x)$ with respect to the scalar product
$\langle\;,\;\rangle_{0,t}$.
Further let us define the dual base of $S_{\lambda}(x;t)$ with
respect to the scalar product $\langle\;,\;\rangle_{q,t}$
by $S_{\lambda}(x;q,t)$ .

We have the following:
\begin{prop}
For a partition $\lambda$, 
we have
\begin{eqnarray}
  |s_{\lambda}\rangle
  &\!\!=\!\!&
  \oint \prod_{j=1}^{\ell(\lambda)}
  \frac{dx_j}{2\pi i x_j}x_j^{\lambda_j} \cdot
  \prod_{i<j}(1-x_i/x_j)\cdot
  \prod_{j=1}^{\ell(\lambda)}
  \exp \left(\sum_{n>0} \frac{a_{-n}}{n}x_j^{-n} \right)|0\rangle,
  \label{Sch1}\\
  |S_{\lambda}(t)\rangle
  &\!\!=\!\!&
  \oint \prod_{j=1}^{\ell(\lambda)}
  \frac{dx_j}{2\pi i x_j}x_j^{\lambda_j} \cdot
  \prod_{i<j}(1-x_i/x_j)\cdot
  \prod_{j=1}^{\ell(\lambda)}
  \exp \left(\sum_{n>0}\frac{1-t^n}{1}\frac{a_{-n}}{n}
  x_j^{-n}\right)|0\rangle,
  \label{Sch2}\\
  |S_{\lambda}(q,t)\rangle
  &\!\!=\!\!&
  \oint \prod_{j=1}^{\ell(\lambda)}
  \frac{dx_j}{2\pi i x_j}x_j^{\lambda_j} \cdot
  \prod_{i<j}(1-x_i/x_j)\cdot
  \prod_{j=1}^{\ell(\lambda)}
  \exp \left(\sum_{n>0}\frac{1}{1-q^n}\frac{a_{-n}}{n}
  x_j^{-n}\right)|0\rangle,
  \label{Sch3}\\
  \langle s_{\lambda}|
  &\!\!=\!\!&
  \oint \prod_{j=1}^{\ell(\lambda)}
  \frac{dx_j}{2\pi i x_j}x_j^{-\lambda_j} \cdot
  \prod_{i<j}(1-x_j/x_i)\cdot
  \langle 0|
  \prod_{j=1}^{\ell(\lambda)}
  \exp \left( \sum_{n>0}\frac{a_{n}}{n}x_j^{n}\right),
  \label{Sch4}\\
  \langle S_{\lambda}(t)|
  &\!\!=\!\!&
  \oint \prod_{j=1}^{\ell(\lambda)}
  \frac{dx_j}{2\pi i x_j}x_j^{-\lambda_j} \cdot
  \prod_{i<j}(1-x_j/x_i)\cdot
  \langle 0|
  \prod_{j=1}^{\ell(\lambda)}
  \exp \left(\sum_{n>0}\frac{1-t^n}{1}\frac{a_{n}}{n}x_j^{n}\right),
  \label{Sch5}\\
  \langle S_{\lambda}(q,t)|
  &\!\!=\!\!&
  \oint \prod_{j=1}^{\ell(\lambda)}
  \frac{dx_j}{2\pi i x_j}x_j^{-\lambda_j} \cdot
  \prod_{i<j}(1-x_j/x_i)\cdot
  \langle 0|
  \prod_{j=1}^{\ell(\lambda)}
  \exp \left( \sum_{n>0}\frac{1}{1-q^n}\frac{a_{n}}{n}x_j^{n}\right).
  \label{Sch6}
\end{eqnarray}
\end{prop}

{\it Proof.}\qquad
An integral representation of the Schur function is well
known\cite{rDJKM}:
\be
  s_{\lambda}(x)=
  \oint \prod_{j=1}^{\ell(\lambda)}
  \frac{dy_j}{2\pi i y_j}y_j^{\lambda_j} \cdot
  \prod_{i<j}(1-y_i/y_j)\cdot
  \prod_{i,j}(1-x_i/y_j)^{-1}.
\ee
Therefore \eq{Sch1} and \eq{Sch4} are correct states.
Since $\langle s_{\lambda}|s_{\mu}\rangle\Bigl|_{t=q}=
\delta_{\lambda,\mu}$, we have
\begin{equation}
  \oint \prod_{j=1}^{\ell(\lambda)}
  \frac{dx_j}{2\pi i x_j} x_j^{-\lambda_j} \cdot
  \prod_{j=1}^{\ell(\mu)}
  \frac{dy_j}{2\pi i y_j} y_j^{\mu_j} \cdot
  \prod_{i<j}(1-x_j/x_i)\cdot
  \prod_{i,j}\frac{1}{1-x_i/y_j}\cdot
  \prod_{i<j}(1-y_i/y_j)
  =\delta_{\lambda,\mu}.
\end{equation}
Note that $S_{\lambda}(x;t)$ is independent of $q$.
Using above identity and proposition 3.5, we can show the following:
\begin{eqnarray*}
  &&
  \langle S_{\lambda}(t),s_{\mu}\rangle_{0,t}=
  \langle S_{\lambda}(t)|s_{\mu}\rangle\Bigl|_{q=0}=
  \delta_{\lambda,\mu},\\
  &&
  \langle S_{\lambda}(q,t),S_{\mu}(t)\rangle_{q,t}=
  \langle S_{\lambda}(q,t)|S_{\mu}(t)\rangle=
  \delta_{\lambda,\mu}. \qquad Q.E.D.
\end{eqnarray*}
We remark that we obtain another expressions of these states by using
another integral representation of the Schur function\cite{rDJKM},
\be
  s_{\lambda'}(x)=(-1)^{|\lambda|}
  \oint \prod_{j=1}^{\ell(\lambda)}
  \frac{dy_j}{2\pi i y_j}y_j^{\lambda_j} \cdot
  \prod_{i<j}(1-y_i/y_j)\cdot
  \prod_{i,j}(1-x_i/y_j).
\ee

Since $K_{\lambda\mu}(q,t)=
\langle S_{\lambda}(q,t),M_{\mu}(q,t)\rangle_{q,t}$,
by using propositions 3.5, 3.6 and 3.9, we can show the following theorem:
\begin{theorem}
Let $\lambda$, $\mu$ be partitions;
$\mu=(\sigma_M^{\rho_M})+\cdots+(\sigma_1^{\rho_1})$.
Kostka matrix $K_{\lambda\mu}(q,t)$ is represented as follows:
\begin{eqnarray}
  &&
  K_{\lambda\mu}(q,t)=\langle S_{\lambda}(q,t)|M_{\mu}(q,t)\rangle \\
  &=\!\!&
  h_{\mu}(q,t) \oint \prod_{j=1}^{\ell(\lambda)}
  \frac{dx_j}{2\pi i x_j}x_j^{-\lambda_j} \cdot
  \prod_{j=1}^{\rho_M} \frac{dy_j}{2\pi i y_j}\cdot
  \prod_{i<j}(1-x_j/x_i)\cdot
  \prod_{i,j} \frac{1}{(x_i/y_j ;q)_{\infty}}\cdot
  F^+_{\mu}(y;q,t) \n
  &=\!\!&
  h_{\mu}(q,t) \oint \prod_{j=1}^{\ell(\lambda)}
  \frac{dx_j}{2\pi i x_j}x_j^{-\lambda_j} \cdot
  \prod_{j=1}^{\rho_M} \frac{dy_j}{2\pi i y_j}\cdot
  \prod_{i<j}(1-x_j/x_i)\cdot
  \prod_{i,j} (-x_i/y_j ;t)_{\infty}\cdot
  F^-_{\mu}(y;q,t).
  \nonumber
\end{eqnarray}
\end{theorem}


\section{Bosonizations of the Integral Formula}
\setcounter{equation}{0}

In the last section, we introduced the Fock space of the boson field
to discuss how to obtain integral formulas for the skew Macdonald
functions.
One may notice, however, the bosonization is something {\it partial}
compared with the case of the Schur functions\cite{rDJKM}
and the Hall-Littlewood functions\cite{rJing},
because we only bosonized the variable $x=x^{N+1}$.
We consider some {\it total} bosonization schemes of the integral
representation formula for the Macdonald symmetric functions
which was obtained in the last section.

\noindent
{\large\bf 4.1}\\
Firstly we treat the case of $\beta\in{\bf Z}_{>0}$. In this case,
we can bosonize the integral formula
by using a similar method to the Jack symmetric function's case
discussed in our previous paper \cite{rAMOS2}.
Let us introduce the following bosonic oscillators having
$A$-type like symmetry:
\be
  [a^a_n, a^b_m] =
  \left\{
  \begin{array}{lll}
     0&{\rm for}& A^{ab}=0, \\
     \displaystyle
     -n \frac{1-t^{|n|}}{1-q^{|n|}}\delta_{n+m,0}&{\rm for}&
     A^{ab}=-1, \\
     \displaystyle
     n \left(\frac{1-t^n}{1-q^n}+\frac{1-t^{-n}}{1-q^{-n}}\right)
     \delta_{n+m,0}&{\rm for}& A^{ab}=2,
  \end{array}\right.
\ee
and $[a^a_0, Q^b  ] = \beta A^{ab}$,
where $n,m\in\bZ$ and $a,b\in\{1,\cdots,N+1\}$.
Here, $A^{ab} = 2\delta^{a,b} - \delta^{a,b+1} - \delta^{a,b-1}$
is the Cartan matrix of $A_{N+1}$ type.
Let us define $A$-type boson fields as follows:
\be
  \phi^a(z)=\phi^a_{\leq 0}(z) + \phi^a_{\geq 0}(z),\quad
  \left\{
  \begin{array}{l}
    \displaystyle
    \phi^a_{\leq 0}(z)=\sum_{n>0}{a^a_{-n} \over n} z^n + Q^a, \\
    \displaystyle
    \phi^a_{\geq 0}(z)=-\sum_{n>0}{a^a_n \over n} z^{-n}+a^a_0 \log z.
  \end{array}
  \right.
\ee
The normal ordering $:~~:$ is defined as moving $\phi_{\geq 0}$ to the
right of $\phi_{\leq 0}$.
%
%
The operator product expansion (in the region $|w/z|<q^{\beta-1}$)
is given as follows:
\be
  \phi^a(z) \phi^b(w)
  \sim
  \left\{\begin{array}{lll}
    0&{\rm for}& A^{ab}=0, \\
    \displaystyle
    \log\biggl(z^{-\beta}\prod_{k=0}^{\beta-1}
    \Bigl(1-q^k w/z\Bigr)^{-1}\biggr)
    &{\rm for}& A^{ab}=-1, \\
    \displaystyle
    \log\biggl((-zw)^\beta q^{-{1\over 2}\beta(\beta-1)}
    \prod_{k=0}^{\beta-1}\Bigl(1-q^k w/z\Bigr)\Bigl(1-q^k z/w\Bigr)
    \biggr)
    &{\rm for}& A^{ab}=2.
  \end{array}\right.
\ee

For $\alpha=(\alpha^1,\cdots,\alpha^{N+1})$,
let $\ket{\alpha}=\exp(\frac{1}{\beta}\sum_{a,b=1}^{N+1}
\alpha^a(A^{-1})_{ab}Q^b)\ket{0}$,
where $\ket{0}$ satisfies $a^a_n\ket{0}=0$
($a=1,\cdots,N+1$ and $n\geq 0$).
This $\ket{\alpha}$ satisfies
$a^a_n|\alpha\rangle=0$ ($a=1,\cdots,N+1$ and $n>0$) and
$a^a_0|\alpha\rangle=\alpha^a|\alpha\rangle$.
We also introduce $\langle\alpha|$ as the dual of $|\alpha\rangle$,
i.e. $\langle\alpha|\alpha'\rangle=\delta_{\alpha,\alpha'}$.

We can state our result as follows:\footnote{
We use the convention
$\prod_{i=1}^n {\cal O}_i={\cal O}_1{\cal O}_2\cdots{\cal O}_n$
for non-commuting ${\cal O}_i$'s.}
\begin{prop}
Let $\lambda$ be defined by {\rm (\ref{deflam})}.
We have the following
$A$-type bosonic realization of the Macdonald symmetric
function for $\beta\in{\bf Z}_{>0}$:
\ba
  P_\lambda(x;q,t)
  &\!\!=\!\!&
  C_{\lambda}^+(q,t)
  \Bigl(-q^{\frac{1}{2}(\beta-1)}\Bigr)^{\frac12
  \beta\sum_{a=1}^Nr_a(r_a-1)} \n
  && \times
  \oint\prod_{a=1}^N \prod_{j=1}^{r_a}
  \frac{dx^a_j}{2\pi i x^a_j }\cdot
  \bra{\tilde\alpha}
  \prod_{a=1}^N \prod_{j=1}^{r_a}
  :\!e^{\phi^a(x^a_j)}\!:\cdot
  \prod_{i=1}^{r_{N+1}}e^{\phi^{N+1}_{\leq 0}(x^{N+1}_i)}
  \ket{\alpha},
\ea
where
$\alpha^a = \beta(r_{a+1}-r_a+1)+s_a$, $(r_{N+2}=0)$,
$\tilde\alpha^a = \alpha^a + \beta\sum_{b=1}^{N+1}A^{ab}r_b$,
and $x_i=x^{N+1}_i$, $r_{N+1}=\infty$ $($ after calculation $)$.
\end{prop}

{\it Proof} \qquad
First we remark that for $\beta\in{\bf Z}_{>0}$ ($t=q^\beta$),
\be
  \Delta(x;q,t)=\prod_{i\neq j}\prod_{k=0}^{\beta-1}
  (1-q^kx_i/x_j),\quad
  \Pi(x,y;q,t)=\prod_{i,j}\prod_{k=0}^{\beta-1}
  \frac{1}{1-q^kx_iy_j}.
\ee
A straightforward calculation of the operator
product expansion shows that this integrand agrees
with that of theorem 3.1. \qquad{\it Q.E.D.}

\noindent
{\large\bf 4.2}\\
Next
we construct another bosonization scheme which is applicable
for the case of $\beta\in{\bf C}$. We utilize Jing's boson field
which was introduced to consider the Hall-Littlewood symmetric
function $P(x;t)$ having one parameter $t$ \cite{rJing}.
Notice that in this case we will not utilize
$A_N$ structure but derive a bosonization formula for $P(x;q,t)$
using finite temperature calculation regarding the parameter $q$ as
playing the role of temperature.
Let us introduce $N$ copy of boson oscillators
$a^a_n$ $(a=1,2,\cdots,N)$,
whose commutation relations are given as follows:
\be
  [a^a_n,a^b_m]=
  n \frac{1}{1-t^{|n|}}\delta_{n+m,0}\delta^{ab}.
\ee
Let ${\cal F}$ be the Fock space of these boson fields:
${\cal F}={\bf Q}(q,t)[a^a_{-1},a^a_{-2},\cdots]\ket{0}$.
Normal ordering $:~~:$ is defined as moving $a^a_n$ to the right
of $a^a_{-n}$ ($n>0$).
We define the grading operator $L_0$ as
\be
  L_0=
  \sum_{a=1}^N \sum_{n=1}^{\infty}(1-t^n)a^a_{-n} a^a_n,
\ee
which satisfies $[L_0,a^a_n]=-na^a_n$.
We introduce boson fields as follows:
\be
  \phi^a(z)=
  -\sum_{n\in{\bf Z}_{\neq 0}}(1-t^{|n|})\frac{a^a_{n}}{n}z^{-n},
  \quad
  \phi^a_-(z)=
  \sum_{n>0}(1-t^n)\frac{a^a_{-n}}{n}z^n.
\ee

Here we state another bosonization formula:
\begin{prop}
Let $\lambda$ be defined by {\rm (\ref{deflam})},
and let $\beta\in{\bf C}$.
We have the following bosonization formula for the
Macdonald symmetric function:
\ba
  P_{\lambda}(x;q,t)
  &\!\!=\!\!&
  C_{\lambda}^+(q,t)(q;q)_{\infty}^{N-\sum_{a=1}^Nr_a}
  (tq;q)_{\infty}^{\sum_{a=1}^Nr_a}
  \prod_{a=1}^Nr_a!\prod_{k=1}^{r_a}\frac{1-t}{1-t^k}
  \label{trace} \\
  &&\times
  \oint\prod_{a=1}^N\prod_{j=1}^{r_a}\frac{dx^a_j}{2\pi i x^a_j}
  \Bigl(x^a_j\Bigr)^{s_a}\cdot
  {\rm Tr}_{{\cal F}}\Biggl( q^{L_0}\;
  \prod_{a=1}^N\prod_{j=1}^{r_a}
  :\!e^{\phi^a(x^a_j)}\!:\cdot
  \prod_{a=1}^N\prod_{j=1}^{r_{a+1}}
  e^{-\phi^a_-(x^{a+1}_j)}
  \Biggr),\nonumber
\ea
where $x_i=x^{N+1}_i$ and $r_{N+1}=\infty$.
\end{prop}

{\it Proof} \qquad
To calculate the trace, we apply
the Clavelli-Shapiro's trace technique\cite{rClSh}.
We introduce the boson oscillators $b^a_n$,
which satisfy the same commutation relation as $a^a_n$ and
commutes with $a^a_m$,
and take the following combinations ($n>0$):
\be
  \tilde{a}^a_n=\frac{a^a_n}{1-q^n}+b^a_{-n},\qquad
  \tilde{a}^a_{-n}=a^a_{-n}+\frac{q^n b^a_n}{1-q^n}.
\ee
Clavelli and Shapiro's argument tells us that
\be
  {\rm Tr}_{\cal F}(q^{L_0} {\cal O})=
  \frac{\bra{0}\tilde{\cal O}\ket{0}}{\prod_{k=1}^{\infty}(1-q^k)^N},
\ee
where ${\cal O}$ is an operator in $a^a_n$, and
$\tilde{{\cal O}}$ is defined as the operator obtained
>from ${\cal O}$ by replacing $a^a_n$ with $ \tilde{a}^a_n$.
Then we obtain
\ba
  :\!e^{\tilde{\phi}^a(z)}\!:
  &\!\!=\!\!&
  \frac{(q;q)_{\infty}}{(qt;q)_{\infty}}
  \exp\left(\sum_{n>0}(1-t^n)\frac{a^a_{-n}}{n}z^n\right)
  \exp\left(-\sum_{n>0}\frac{1-t^n}{1-q^n}
  \frac{a^a_n}{n}z^{-n}\right) \n
  &&\qquad\quad\times
  \exp\left(-\sum_{n>0}(1-t^n)\frac{b^a_{-n}}{n}z^{-n}\right)
  \exp\left(\sum_{n>0}\frac{(1-t^n)q^n}{1-q^n}
  \frac{b^a_n}{n}z^n \right),
\ea
and
\be
  :\!e^{-\tilde{\phi}^a_-(z)}\!:
  \;=
  \exp\left(-\sum_{n>0}(1-t^n)\frac{a^a_{-n}}{n}z^n\right)
  \exp\left(-\sum_{n>0}\frac{(1-t^n)q^n}{1-q^n}
  \frac{b^a_n}{n}z^n\right).
\ee
We have the following OPE's;
in the region $q<|z_2/z_1|<1$,
\be
  :\!e^{\tilde{\phi}^a(z_1)}\!:\::\!e^{\tilde{\phi}^b(z_2)}\!:
  \;=
  \left\{\begin{array}{lll}
    \displaystyle
    :\!e^{\tilde{\phi}^a(z_1)+\tilde{\phi}^a(z_2)}\!:
    \frac{(z_2/z_1;q)_{\infty}}{(t z_2/z_1;q)_{\infty}}
    \frac{(q z_1/z_2;q)_{\infty}}{(qt z_1/z_2;q)_{\infty}}&
    {\rm for} & a=b, \\
    :\!e^{\tilde{\phi}^a(z_1)+\tilde{\phi}^b(z_2)}\!:&
    {\rm for} & a\neq b,
  \end{array}\right.
\ee
and in the region $|w/z|<1$,
\be
  :\!e^{\tilde{\phi}^a(z)}\!:\::\!e^{-\tilde{\phi}^b_-(w)}\!:
  \;=
  \left\{\begin{array}{lll}
    \displaystyle
    :\!e^{\tilde{\phi}^a(z)-\tilde{\phi}^a_-(w)}\!:
    \frac{(tw/z;q)_{\infty}}{(w/z;q)_{\infty}}&
    {\rm for} & a=b, \\
    :\!e^{\tilde{\phi}^a(z)-\tilde{\phi}^b_-(w)}\!:&
    {\rm for} & a\neq b.
  \end{array}\right.
\ee

By using these equations,
$(q;q)_{\infty}^{N-\sum_{a=1}^Nr_a}
(tq;q)_{\infty}^{\sum_{a=1}^Nr_a}
{\rm Tr}_{{\cal F}}(\cdots)$
in \eq{trace} becomes
\ba
  &&
  \prod_{a=1}^N\prod_{i,j=1 \atop i<j}^{r_a}
  \frac{(x^a_j/x^a_i;q)_{\infty}}{(tx^a_j/x^a_i;q)_{\infty}}
  \frac{(q x^a_i/x^a_j;q)_{\infty}}{(qt x^a_i/x^a_j;q)_{\infty}}
  \cdot
  \prod_{a=1}^N\prod_{i=1}^{r_a}\prod_{j=1}^{r_{a+1}}
  \frac{(tx^{a+1}_i/x^a_j;q)_{\infty}}{(x^{a+1}_i/x^a_j;q)_{\infty}}\n
  &=\!\!&
  \prod_{a=1}^N\Delta(x^a;q,t)\Pi(x^{a+1},\overline{x^a};q,t)\cdot
  \prod_{a=1}^N\prod_{i,j=1 \atop i<j}^{r_a}
  \frac{(1-tx^a_i/x^a_j)}{(1-x^a_i/x^a_j)}.
\ea
For each $a$ and a permutation $\sigma$,
we change the integration variables $x^a_i\rightarrow x^a_{\sigma(i)}$.
Then by using the identity\cite{rMac}
\be
  \sum_{\sigma\in S_n}\prod_{1\leq i<j\leq n}
  \frac{x_{\sigma(i)}-tx_{\sigma(j)}}
       {x_{\sigma(i)}- x_{\sigma(j)}}
  =
  \prod_{k=1}^n \frac{1-t^k}{1-t},
\ee
where $S_n$ is the $n$-th symmetric group,
the integrand agrees with that of theorem 3.1. \qquad{\it Q.E.D.}

\section{Discussion}
\setcounter{equation}{0}

In this paper we have obtained integral representations of the
(skew-)Macdonald symmetric functions (theorems 3.2, 3.4 and corollary 3.8)
and their boson realizations (propositions 4.1, 4.2 and theorem 3.7).
The two maps in the proposition 3.1 have played
an essential role in our derivation.

Our first physical motivation for this study is
calculation of the correlation functions of the Calogero-Sutherland
model. The results obtained in this paper
and ref.\ \cite{rAMOS1,rMY1,rMY2,rAMOS2}
will help us to do it. 
In particular skew Jack symmetric functions will be useful for higher
point correlation functions.
Of course, concerning the analysis for the Calogero-Sutherland model,
the Macdonald symmetric functions are unnecessary, but
sometimes calculation for $q$-deformed quantities is
more transparent than the original ones.
We have also constructed free boson realizations for the integral
representations. These realizations will also help us
in calculation for correlation functions.
However, in comparison with the case of the Jack symmetric functions,
these free field expressions are ad hoc in the sense that
they merely give the desired integrands of the integral
representations (see the next paragraph).
Another motivation is to solve the Ruijsenaars model, i.e.
model with elliptic potential.
At present this problem seems to be difficult yet.

For mathematical interest, we would like to mention
the relation between free field realizations and symmetry algebras.
In the case of the Jack symmetric function\cite{rAMOS1,rMY1,rAMOS1,rAMOS2},
this function is realized on the boson Fock space
as the state obtained by the action of screening currents
$:\!e^{\alpha_\pm \phi^a(z)}\!:$ on the vacuum.
This state is the singular vector of the $W_N$ algebra.
On the other hand, in the free boson realization,
the $W_N$ algebra is the commutant of these screening currents
$:\!e^{\alpha_\pm \phi^a(z)}\!:$.
So we have the following natural question;
in the case of the Macdonald symmetric functions,
what algebra appears as the commutant of the vertex operators
$:\!e^{\phi^a(x)}\!:$ used in section 4 ?

After finishing this work, we knew that
Frenkel and Reshetikhin constructed certain $q$-deformations of
the Virasoro and $W$-algebras \cite{rFR} by
utilizing the free boson realization of the
quantum affine algebra $U_q(\widehat{sl}_N)$ studied in ref.\ \cite{rAOS}.
It seems interesting to clarify the connection
between our vertex operators introduced in section 4
and the $q$-deformed algebras.

\vskip 5mm
\noindent{\bf Acknowledgments:}

We would like to thank Y.~Matsuo, K.~Iohara, M.~Noumi and Y.~Yamada
for valuable discussions.
J.S. is very grateful
to T. Eguchi
for discussions and encouragements, and
to A.~LeClair and V.~Pasquier
for their kind hospitality at Cornell Univ. and Saclay.
This work is supported in part by Grant--in--Aid for Scientific
Research from Ministry of Science and Culture.


\end{document}